# A Projection-Based Reformulation and Decomposition Algorithm for Global Optimization of a Class of Mixed Integer Bilevel Linear Programs

Dajun Yue · Jiyao Gao · Bo Zeng · Fengqi You[*]

**Abstract** We propose an extended variant of the reformulation and decomposition algorithm for solving a special class of mixed-integer bilevel linear programs (MIBLPs) where continuous and integer variables are involved in both upper- and lower-level problems. In particular, we consider MIBLPs with upper-level constraints that involve lower-level variables. We assume that the inducible region is nonempty and all variables are bounded. By using the reformulation and decomposition scheme, an MIBLP is first converted into its equivalent single-level formulation, then computed by a column-and-constraint generation based decomposition algorithm. The solution procedure is enhanced by a projection strategy that does not require the relatively complete response property. To ensure its performance, we prove that our new method converges to the global optimal solution in a finite number of iterations. A large-scale computational study on random instances and instances of hierarchical supply chain planning are presented to demonstrate the effectiveness of the algorithm.

D. Yue
Northwestern University, Evanston, Illinois 60208, USA.
J. Gao
Cornell University, Ithaca, New York 14853, USA.
B. Zeng
University of Pittsburgh, Pittsburgh, Pennsylvania 15261, USA.
F. You (corresponding author)
Cornell University, Ithaca, New York 14853, USA. Tel: +1 607 882 5530, e-mail: fengqi.you@cornell.edu

# 1. Introduction

We present an algorithm for solving a special class of mixed-integer bilevel linear programs (MIBLPs) of the following form.

$$\textbf{(P0)} \quad \min_{x^u, y^u, x^{l0}, y^{l0}} \quad c_R^t x^u + c_Z^t y^u + d_R^t x^{l0} + d_Z^t y^{l0} \tag{1}$$

$$\text{s.t.} \quad A_R x^u + A_Z y^u + B_R x^{l0} + B_Z y^{l0} \le r \tag{2}$$

$$x^u \in \mathbb{R}_+^{m_R}, y^u \in \mathbb{Z}_+^{m_Z}, x^{l0} \in \mathbb{R}_+^{n_R}, y^{l0} \in \mathbb{Z}_+^{n_Z} \tag{3}$$

$$\left(x^{l0}, y^{l0}\right) \in \underset{\left(x^l, y^l\right)}{\text{argmax}} \left\{ \quad w_R^t x^l + w_Z^t y^l : \right. \tag{4}$$

$$P_R x^l + P_Z y^l \le s - Q_R x^u - Q_Z y^u \tag{5}$$

$$\left. x^l \in \mathbb{R}_+^{n_R}, y^l \in \mathbb{Z}_+^{n_Z} \quad \right\} \tag{6}$$

In formulation (P0), there is an ambiguity when multiple lower-level optimal solutions exist [1]. In the optimistic (or strong) formulation, the lower-level decision maker selects $\left(x^{l0}, y^{l0}\right)$ from his optimal solution set according to the interests of the upper-level decision maker [2]. On the contrary, in the pessimistic (or weak) formulation, the lower-level decision maker select one optimal solution to against the upper level decision maker's interest [3]. In this paper, the optimistic formulation is treated. Without loss of generality, the lower-level program in (P0) can be converted into a minimization problem by changing the sign of the lower-level objective function.

Bilevel programs [4], including MIBLPs, are frequently utilized to model Stackelberg games in game theory [5,6]. MIBLPs are intrinsically challenging to solve, and the use of mixed-integer linear programming (MILP) algorithms for solving MIBLPs is not straightforward [7-9]. Although MIBLPs can be solved using some general-purpose mixed-integer bilevel nonlinear program (MIBNLP) algorithms, most of them only have small-scale applications reported in the literature [10-12]. For example, the global MIBNLP algorithm proposed by Mitsos [13] was applied by Fliscounakis et al. to solve a specific large-scale MIBLP problem for a power system [14]. However, the detailed model formulation and the resulting problem size of this MIBLP are not reported. Most existing MIBLP algorithms are proposed to handle special classes of (P0), such as integer bilevel linear programs [15], MIBLPs with special constraint structures [16], MIBLPs without continuous upper-level variables [7,17,18], and/or MIBLPs without continuous lower-level variables [19,20]. Fischetti et al. [21,22] introduced a new general-purpose algorithm for MIBLPs based on a branch-and-cut framework, where new classes of valid inequalities and effective preprocessing procedures are introduced. We mention that Zeng and An [23] proposed the reformulation and decomposition method to solve MIBLPs with continuous and integer variables in both upper- and lower-level programs, which provides a rather general strategy and framework to attack those difficult problems. In a capacity expansion planning problem [24], this new method demonstrates a very strong solution capacity. Nevertheless, the original reformulation and decomposition method did not consider MIBLPs with upper-level constraints involving lower-level variables (called connecting constraints [25] in what follows).

In this work, we extend the original reformulation and decomposition scheme in [23] to solve MIBLPs in a more general form, i.e., the form of (P0) that has connecting constraints. Specifically, the decomposition algorithm, i.e., the master and subproblems are modified or updated according to the structure of (P0). Also, the finite convergence proof is generalized to guarantee that our new development derives (P0)'s global optimal solutions. We point out that a new feature based on a novel projection-based formulation is introduced to handle the case where the relatively completely response property does not hold. This enhancement could be of a critical value as many real problems may not have that property. To verify this new development of the reformulation and decomposition method, a large-scale computational study on two types of random instances and instances of hierarchical supply chain planning are presented. Similarly to [24], our algorithm demonstrates a desirable computational capacity on more general MIBLP's.

The rest of the paper is organized as follows. A brief literature review is given in Section 2. We provide the preliminaries in Section 3. We describe the reformulation procedure, especially the projection-based formulation, in Section 4. The decomposition algorithm is presented in Section 5. We discuss several implementation issues in Section 6. A rather comprehensive computational study are presented in Section 7. We conclude the article in Section 8.

## 2. Literature review

A variety of approaches have been proposed to solve MIBLP problems in the literature. Moore and Bard [8,26] proposed the first branch-and-bound algorithms for MIBLPs. Dempe [19] and Hemmati and Smith [16] proposed a cutting plane approach. Saharidis and Ierapetritou [27] proposed an algorithm based on Benders decomposition. DeNegre and Ralphs [15] presented a branch-and-cut algorithm. Köppe et al. [20] proposed a parametric integer programming algorithm. Recently, Xu and Wang [7] developed an exact algorithm based on the branch-and-bound framework. Fischetti et al. [21,22] introduced new classes of linear inequalities in a branch-and-cut framework. Poirion et al. [28] proposed a cut-generation algorithm and row-and-column generation framework. As mentioned, Zeng and An [23] proposed the original reformulation-and-decomposition scheme. These algorithms have been proposed to handle different classes of MIBLP problems. This work contributes to solving a class of MIBLPs in the form of (P0).

The relevant literature also includes studies on bilevel nonlinear programs (BNLP) and MIBNLP. Edmunds and Bard [29] proposed a branch-and-bound algorithm for MIBNLPs. Gümüş and Floudas [10] proposed a vertex polyhedral convex hull representation. Faísca et al. [30] and Domínguez and Pistikopoulos [11] employed parametric programming approaches. Mitsos et al. [31] and Mitsos [13] proposed bounding algorithms for global optimization of BNLP and MIBNLP problems. Kleniati and Adjiman [32,12] proposed branch-and-sandwich algorithms for solving BNLP and MIBNLP problems. Based on the computational performances reported in the literature mentioned above, the computational performance of existing MIBNLP algorithms in solving medium to large MIBLP problems needs to be further tested.

Algorithms have also been proposed for solving other types of relevant programs, including min-max programs [33,34,9], semi-infinite programs [35-37], and generalized semi-infinite programs [38-40]. However, there is no direct relationship between the proposed algorithm and these algorithms. It is worth mentioning that there is a large body of literature on heuristic and meta-heuristic algorithms for bilevel optimization problems [41,42], which is out of the scope of this work.

## 3. Preliminaries

In this section, we introduce some definitions and an assumption that will be used in the algorithm.

**Definition 1** We denote $\Omega$:

$$\Omega = \left\{ \left( x^u, y^u, x^{l0}, y^{l0} \right) : \quad \begin{array}{l} A_R x^u + A_Z y^u + B_R x^{l0} + B_Z y^{l0} \leq r, \\ Q_R x^u + Q_Z y^u + P_R x^{l0} + P_Z y^{l0} \leq s, \\ x^u \in \mathbb{R}_+^{m_R}, y^u \in \mathbb{Z}_+^{m_Z}, x^{l0} \in \mathbb{R}_+^{n_R}, y^{l0} \in \mathbb{Z}_+^{n_Z} \end{array} \right\} \tag{7}$$

the *MIBLP constraint region*.

**Definition 2** For any given $\left( x^u, y^u \right) \in \mathbb{R}_+^{m_R} \times \mathbb{Z}_+^{m_Z}$, we denote $\Omega_{\left( x^u, y^u \right)}$:

$$\Omega_{\left(x^u, y^u\right)} = \left\{ \left(x^{l0}, y^{l0}\right) : P_R x^{l0} + P_Z y^{l0} \leq s - Q_R x^u - Q_Z y^u, x^{l0} \in \mathbb{R}_+^{n_R}, y^{l0} \in \mathbb{Z}_+^{n_Z} \right\} \tag{8}$$

the *lower-level feasible region*.

**Definition 3** For any given $\left(x^u, y^u\right) \in \mathbb{R}_+^{m_R} \times \mathbb{Z}_+^{m_Z}$, we denote $M_{\left(x^u, y^u\right)}$:

$$M_{\left(x^u, y^u\right)} = \arg\max_{\left(x^{l0}, y^{l0}\right)} \left\{ w_R^t x^{l0} + w_Z^t y^{l0} : \left(x^{l0}, y^{l0}\right) \in \Omega_{\left(x^u, y^u\right)} \right\} \tag{9}$$

the *lower-level rational reaction set*.

**Definition 4** We denote *IR*:

$$IR = \left\{ \left(x^u, y^u, x^{l0}, y^{l0}\right) : \left(x^u, y^u, x^{l0}, y^{l0}\right) \in \Omega, \ \left(x^{l0}, y^{l0}\right) \in M_{\left(x^u, y^u\right)} \right\} \tag{10}$$

the *inducible region*, which represents the feasible region at the upper-level program.

Hence, a more general definition for (P0) can be:

**(P1)** $\quad \min\limits_{x^u, y^u, x^{l0}, y^{l0}} \quad c_R^t x^u + c_Z^t y^u + d_R^t x^{l0} + d_Z^t y^{l0} : \left(x^u, y^u, x^{l0}, y^{l0}\right) \in IR$

**Assumption 1** The inducible region *IR* is nonempty and all variables have finite bounds.

Assumption 1 ensures that the feasible set of MIBLP problem (P0) or (P1) is nonempty. We note that the connecting constraints (2) in the upper level problem can make the inducible region empty even if the lower level problem has an optimal solution for the selection in the upper level problem. In cases that an optimal solution does not exist due to open inducible region or non-lower-semi-continuous optimal value function [43-45], we are interested in the infimum of the objective function and $\epsilon$-optimal solutions as shown in [23]. In this work, we do not explicitly consider the cases of unboundedness and infeasibility.

## 4. Reformulations

### 4.1. Optimal value reformulation

First, we follow the convention to reformulate the optimistic MIBLP (P0) using the optimal value transformation [46], as given below.

**(P2)** $\min\limits_{x^u, y^u, x^{l0}, y^{l0}} \quad c_R^t x^u + c_Z^t y^u + d_R^t x^{l0} + d_Z^t y^{l0}$ $\hfill$ (11)

s.t. $A_R x^u + A_Z y^u + B_R x^{l0} + B_Z y^{l0} \leq r$ $\hfill$ (12)

$\qquad Q_R x^u + Q_Z y^u + P_R x^{l0} + P_Z y^{l0} \leq s$ $\hfill$ (13)

$$w_R^l x^{l0} + w_Z^l y^{l0} \geq \max_{x^l, y^l} \left\{ \begin{array}{l} w_R^l x^l + w_Z^l y^l : \\ Q_R x^u + Q_Z y^u + P_R x^l + P_Z y^l \leq s \\ x^l \in \mathbb{R}_+^{n_R}, y^l \in \mathbb{Z}_+^{n_Z} \end{array} \right\} \tag{14}$$

$$x^u \in \mathbb{R}_+^{m_R}, y^u \in \mathbb{Z}_+^{m_Z}, x^{l0} \in \mathbb{R}_+^{n_R}, y^{l0} \in \mathbb{Z}_+^{n_Z} \tag{15}$$

Because the lower-level program is a maximization problem, constraint (14) ensures that $\left( x^{l0}, y^{l0} \right)$ is an optimal solution to the lower-level problem, for any given $\left( x^u, y^u \right)$. Furthermore, (P2) corresponds to the optimistic formulation because $\left( x^{l0}, y^{l0} \right)$ are controlled for the benefits of the upper-level program.

### 4.2. Projection-based single-level formulation

It is a standard approach to reduce bilevel programs to equivalent single-level programs [47-50]. In the case that the lower-level program is a linear program (LP), one can replace the lower-level program with its corresponding Karush-Kuhn-Tucker (KKT) conditions. However, the discrete variables $y^l$ in (P2) renders the lower level problem non-convex and hence cannot be replaced by its KKT conditions. To deeply analyze their impact, we introduce the following projection concept.

**Definition 5** We denote $\text{Proj}_{y^l} \Omega$ :

$$\text{Proj}_{y^l} \Omega = \left\{ y^l \in \mathbb{Z}_+^{n_Z} : \exists \left( x^u, y^u, x^l \right) \in \mathbb{R}_+^{m_R} \times \mathbb{Z}_+^{m_Z} \times \mathbb{R}_+^{n_R} \text{ with } \left( x^u, y^u, x^l, y^l \right) \in \Omega \right\} \tag{16}$$

the *projection of the constraint region on the space of lower-level integer variables*, which represents the collection of all admissible $y^l$.

Following the idea in [23], we can separate the continuous and integer variables in the lower-level program and restructure the right-hand-side of (14):

$$w_R^l x^{l0} + w_Z^l y^{l0} \geq \max_{y^l \in \text{Proj}_{y^l} \Omega} \left( w_Z^l y^l + \max_{x^l} \left\{ w_R^l x^l : P_R x^l \leq s - P_Z y^l - Q_R x^u - Q_Z y^u, x^l \in \mathbb{R}_+^{n_R} \right\} \right) \tag{17}$$

As pointed in [23], because the second maximization problem in (17) is an LP, we can replace it with the KKT-conditions, thus having the following equivalent form:

$$w_R^t x^{l0} + w_Z^t y^{l0} \geq \max_{y^l \in \text{Proj}_{y^l}\Omega} w_Z^t y^l + w_R^t x^l \tag{18}$$

$$\text{s.t. } \left(x^l, \pi\right) \in \begin{cases} P_R x^l \leq s - Q_R x^u - Q_Z y^u - P_Z y^l \\ P_R^t \pi \geq w_R^t, x^l \perp \left(P_R^t \pi - w_R^t\right) \\ \pi \perp \left(s - Q_R x^u - Q_Z y^u - P_R x^l - P_Z y^l\right) \\ x^l \in \mathbb{R}_+^{n_R}, \pi \in \mathbb{R}_+^{n_L} \end{cases} \tag{19}$$

where the $\perp$ (perpendicular) operator enforces the perpendicularity condition between the vectors on the left- and right-hand sides, i.e., their element-by-element product is equal to zero. According to Assumption 1, all lower-level integer variables are bounded. Hence, $\text{Proj}_{y^l}\Omega$ is a finite set. For ease of exposition, we let $Y^L = \left\{y^{l,1}, y^{l,2}, ..., y^{l,J}\right\}$ (indexed by $j$) represents the finite set of all $y^l$ such that $y^l \in \text{Proj}_{y^l}\Omega$.

Then, following the strategy of [23], by enumerating $y^{l,j} \in Y^L$ and introducing corresponding primal and dual variables $\left(x^{l,j}, \pi^j\right)$ and their related KKT-conditions, we have the next formulation (P3). Though finite, the set $Y^L$ could be extremely large and cannot be bounded by a polynomial in the dimension of the problem.

$$\textbf{(P3)} \min_{\substack{x^u, y^u, x^{l0}, y^{l0} \\ x^{l,j}, \pi^{l,j}}} c_R^t x^u + c_Z^t y^u + d_R^t x^{l0} + d_Z^t y^{l0} \tag{11}$$

$$\text{s.t. } \quad (12), (13), \text{ and } (15)$$

$$w_R^t x^{l0} + w_Z^t y^{l0} \geq w_R^t x^{l,j} + w_Z^t y^{l,j}, \quad \forall y^{l,j} \in Y^L \tag{20}$$

$$P_R x^{l,j} \leq s - Q_R x^u - Q_Z y^u - P_Z y^{l,j}, \quad \forall y^{l,j} \in Y^L \tag{21}$$

$$P_R^t \pi^j \geq w_R^t, x^{l,j} \perp \left(P_R^t \pi^j - w_R^t\right), \quad \forall y^{l,j} \in Y^L \tag{22}$$

$$\pi^j \perp \left(s - Q_R x^u - Q_Z y^u - P_R x^{l,j} - P_Z y^{l,j}\right), \quad \forall y^{l,j} \in Y^L \tag{23}$$

$$x^{l,j} \in \mathbb{R}_+^{n_R}, \pi^j \in \mathbb{R}_+^{n_L}, \quad \forall y^{l,j} \in Y^L \tag{24}$$

We mention that it may seem that (P3) is equivalent to (P2). However, this is not always true. As epitomized by the following example, (P3) can be infeasible even when (P2) is feasible.

$$\textbf{(Q0)} \min_{y^u, x^{l0}, y^{l0}} \quad -y^u - y^{l0} \tag{25}$$

$$\text{s.t. } \quad y^u \in \{0,1\}, x^{l0} \in \mathbb{R}_+, y^{l0} \in \{0,1\} \tag{26}$$

$$\left(x^{l0}, y^{l0}\right) \in \underset{x^l, y^l}{\arg\max} \quad -x^l - y^l \tag{27}$$

$$\text{s.t. } \quad x^l - y^l \leq -y^u \tag{28}$$

$$x^l \in \mathbb{R}_+, y^l \in \{0,1\} \tag{29}$$

It is easy to see that the upper-level optimal solution is $\left(y^{u*}, x^{l0*}, y^{l0*}\right) = (1,0,1)$ with an optimal value of $-2$. Now, let $y^{l,1} = 1$ and $y^{l,2} = 0$. If we reformulate (Q0) according to the formulation in (P3), we then have:

(**Q1**) $\min\limits_{\substack{y^u, x^{l0}, y^{l0} \\ x^{l,1}, \pi^1, x^{l,2}, \pi^2}} - y^u - y^{l0}$ (25)

s.t. (26)

$$y^u + x^{l0} - y^{l0} \le 0 \tag{30}$$

$$-x^{l0} - y^{l0} \ge -x^{l,1} - 1 \tag{31}$$

$$x^{l,1} \le -y^u + 1 \tag{32}$$

$$\pi^1 \ge -1, \quad x^{l,1} \perp \left(\pi^1 + 1\right) \tag{33}$$

$$\pi^1 \perp \left(-y^u - x^{l,1} + 1\right) \tag{34}$$

$$x^{l,1}, \pi^1 \in \mathbb{R}_+ \tag{35}$$

$$-x^{l0} - y^{l0} \ge -x^{l,2} \tag{36}$$

$$x^{l,2} \le -y^u \tag{37}$$

$$\pi^2 \ge -1, \quad x^{l,2} \perp \left(\pi^2 + 1\right) \tag{38}$$

$$\pi^2 \perp \left(-y^u - x^{l,2}\right) \tag{39}$$

$$x^{l,2}, \pi^2 \in \mathbb{R}_+ \tag{40}$$

where (31) – (35) correspond to $y^{l,1} = 1$ and (36) – (39) correspond to $y^{l,2} = 0$. It is interesting to note that $\left(y^{u*}, x^{l0*}, y^{l0*}\right) = (1,0,1)$ is infeasible to (Q1) because constraint (37) indicates $x^{l,2} \le -1$, which contradicts with $x^{l,2} \in \mathbb{R}_+$ in (40). Hence, formulation (P3) is not equivalent to (P2). As noted in [23], this issue is caused by the lack of relatively complete response property in (Q0).

**Definition 6** We denote $\text{Proj}_{\left(x^u, y^u\right)} \Omega$:

$$\text{Proj}_{\left(x^u, y^u\right)} \Omega = \left\{ \left(x^u, y^u\right) \in \mathbb{R}_+^{m_R} \times \mathbb{Z}_+^{m_Z} : \exists \left(x^{l0}, y^{l0}\right) \in \mathbb{R}_+^{n_R} \times \mathbb{Z}_+^{n_Z} \text{ with } \left(x^u, y^u, x^{l0}, y^{l0}\right) \in \Omega \right\} \tag{41}$$

the *projection of the constraint region on the space of upper-level variables*.

Next, we re-define the relatively complete response property in [23] using our projection concepts.

**Definition 7** An MIBLP in the form of (P0) has the *relatively complete response* property if for any 3-tuple $\left( x^u, y^u, y^l \right)$ such that $\left( x^u, y^u \right) \in \text{Proj}_{\left( x^u, y^u \right)} \Omega$ and $y^l \in \text{Proj}_{y^l} \Omega$, the following LP is feasible and has a finite optimal value.

$$\max_{x^l} \quad w_R^t x^l \tag{42}$$

$$\text{s.t.} \quad P_R x^l \leq s - Q_R x^u - Q_Z y^u - P_Z y^l \tag{43}$$

$$x^l \in \mathbb{R}_+^{n_R} \tag{44}$$

Taking (Q0) as example, since there is no upper-level continuous variable, we have $\text{Proj}_{y^u} \Omega = \{0,1\}$ and $\text{Proj}_{y^l} \Omega = \{0,1\}$. When the couple $\left( y^u, y^l \right)$ is fixed at $(1,1)$, $(0,1)$ or $(0,0)$, the related LP $(27) - (29)$ is feasible and has a finite optimal value. However, when the couple $\left( y^u, y^l \right)$ is fixed at $(1,0)$ the related LP $(27) - (29)$ becomes infeasible. Hence, problem (Q0) does not have the relatively complete response property.

In formulation (Q1), constraints $(36) - (40)$ corresponding to $y^{l,2} = 0$ are imposed regardless of the value of $y^u$, which excludes any 3-tuple $\left( y^u, x^{l0}, y^{l0} \right)$ such that $y^u = 1$. This is incorrect according to the previous analysis. A straightforward approach to fix this issue is not to impose constraints $(36) - (40)$ when $y^u = 1$, thus having

**(Q2)** $\min\limits_{\substack{y^u, x^{l0}, y^{l0} \\ x^{l,1}, \pi^1, x^{l,2}, \pi^2}} - y^u - y^{l0}$ \hfill (25)

$$\text{s.t.} \quad (26) \text{ and } (30)$$

$$\left[ y^u \in \{0,1\} \right] \Rightarrow [ (31) - (35) ] \tag{45}$$

$$\left[ y^u = 0 \right] \Rightarrow [ (36) - (40) ] \tag{46}$$

Next, we generalize the equivalent reformulation from (Q0) to (Q2).

**Definition 8** For any given $y^{l,j} \in Y^L$, we denote $P\left( y^{l,j} \right)$

$$P\left( y^{l,j} \right) = \left\{ \begin{array}{c} \left( x^u, y^u, x^l \right) : Q_R x^u + Q_Z y^u + P_R x^l \leq s - P_Z y^{l,j}, \\ x^u \in \mathbb{R}_+^{m_R}, y^u \in \mathbb{Z}_+^{m_Z}, x^l \in \mathbb{R}_+^{n_R} \end{array} \right\} \tag{47}$$

**Definition 9** For any given $y^{l,j} \in Y^L$, we denote $\text{Proj}_{\left( x^u, y^u \right)} P\left( y^{l,j} \right)$

$$\text{Proj}_{\left( x^u, y^u \right)} P\left( y^{l,j} \right) = \left\{ \left( x^u, y^u \right) \in \mathbb{R}_+^{m_R} \times \mathbb{Z}_+^{m_Z} : \exists x^l \in \mathbb{R}_+^{n_R} \text{ with } \left( x^u, y^u, x^l \right) \in P\left( y^{l,j} \right) \right\} \tag{48}$$

With Definitions 8 and 9, we propose a projection-based formulation (P4), which is equivalent to (P0), even when (P0) does not have the relatively complete response property. Hence, it extends the scope investigated in [23] to a more general situation, which actually has an important value as many real problems may not have that property. Notably, although the formulation of (P4) may seem similar to the subroutine 1 adopted in [13], the concepts and approaches of these two ideas are essentially different. The subroutine 1 in [13] is iteratively deriving a bound box on upper level variables based on a complete lower level solution, while (P4) is built on the projection of a feasible discrete variable $y^{l,j}$ on upper level variables.

$$\textbf{(P4)} \min_{\substack{x^u, y^u, x^{l0}, y^{l0} \\ x^{l,j}, \pi^{l,j}}} c_R^t x^u + c_Z^t y^u + d_R^t x^{l0} + d_Z^t y^{l0} \tag{11}$$

$$\text{s.t.} \quad (12), (13), \text{ and } (15)$$

$$\left[\left(x^u, y^u\right) \in \text{Proj}_{\left(x^u, y^u\right)} P\left(y^{l,j}\right)\right]$$

$$\Rightarrow \begin{bmatrix} w_R^t x^u + w_Z^t y^{l0} \geq w_R^t x^{l,j} + w_Z^t y^{l,j} \\ P_R x^{l,j} \leq s - Q_R x^u - Q_Z y^u - P_Z y^{l,j} \\ P_R^t \pi^j \geq w_R^t, x^{l,j} \perp \left(P_R^t \pi^j - w_R^t\right) \\ \pi^j \perp \left(s - Q_R x^u - Q_Z y^u - P_R x^{l,j} - P_Z y^{l,j}\right) \\ x^{l,j} \in \mathbb{R}_+^{n_R}, \pi^j \in \mathbb{R}_+^{n_L} \end{bmatrix}, \quad \forall y^{l,j} \in Y^L \tag{49}$$

In formulation (P4), constraint (49) indicates that constraints (20) – (24) corresponding to any given $y^{l,j} \in Y^L$ will only be imposed if $\left(x^u, y^u\right) \in \text{Proj}_{\left(x^u, y^u\right)} P\left(y^{l,j}\right)$. In the following, a formal argument is presented to show equivalence between our projection-based formulation and (P0).

**Theorem 1** *The projection-based single-level formulation* (P4) *is equivalent to the original MIBLP problem* (P0).

*Proof* Since it has been shown that (P0) is equivalent to (P2) [23,31], we now prove (P4) is also equivalent to (P2). It is sufficient to show constraint (49) is equivalent to (14) since the other constraints and the objectives are the same in both problems.

Let $\left(\dot{x}^u, \dot{y}^u, \dot{x}^{l0}, \dot{y}^{l0}\right)$ be a feasible solution to (P2), we then have

$$w_R^t \dot{x}^{l0} + w_Z^t \dot{y}^{l0} \geq \max_{x^l, y^l} \left\{ w_R^t x^l + w_Z^t y^l : P_R x^l + P_Z y^l \leq s - Q_R \dot{x}^u - Q_Z \dot{y}^u, x^l \in \mathbb{R}_+^{n_R}, y^l \in \mathbb{Z}_+^{n_Z} \right\} \tag{50}$$

$$\Leftrightarrow \begin{bmatrix} w_R^t \dot{x}^{l0} + w_Z^t \dot{y}^{l0} \geq \max w_R^t x^{l,j} + w_Z^t y^{l,j} : \\ \text{s.t. } P_R x^{l,j} \leq s - Q_R \dot{x}^u - Q_Z \dot{y}^u - P_Z y^{l,j}, \\ x^{l,j} \in \mathbb{R}_+^{n_R} \end{bmatrix},$$
$$\forall y^{l,j} \in \left\{ \begin{array}{l} y^l \in \mathbb{Z}_+^{n_Z} : \\ \exists x^l \in \mathbb{R}_+^{n_R} \text{ with } \left( x^l, y^l \right) \in \Omega_{\left( \dot{x}^u, \dot{y}^u \right)} \end{array} \right\} \tag{51}$$

$$\Leftrightarrow \begin{bmatrix} w_R^t \dot{x}^{l0} + w_Z^t \dot{y}^{l0} \geq w_R^t x^{l,j} + w_Z^t y^{l,j} \\ P_R x^{l,j} \leq s - Q_R \dot{x}^u - Q_Z \dot{y}^u - P_Z y^{l,j} \\ P_R^t \pi^j \geq w_R, x^{l,j} \perp \left( P_R^t \pi^j - w_R^t \right) \\ \pi^j \perp \left( s - Q_R \dot{x}^u - Q_Z \dot{y}^u - P_R x^{l,j} - P_Z y^{l,j} \right) \\ x^{l,j} \in \mathbb{R}_+^{n_R}, \pi^j \in \mathbb{R}_+^{n_L} \end{bmatrix}, \quad \forall y^{l,j} \in \left\{ \begin{array}{l} y^l \in \mathbb{Z}_+^{n_Z} : \\ \exists x^l \in \mathbb{R}_+^{n_R} \text{ with } \left( x^l, y^l \right) \in \Omega_{\left( \dot{x}^u, \dot{y}^u \right)} \end{array} \right\} \tag{52}$$

Only the constraints (49) corresponding to $y^{l,j} \in \left\{ y^l \in \mathbb{Z}_+^{n_Z} : \exists x^l \in \mathbb{R}_+^{n_R} \text{ with } \left( x^l, y^l \right) \in \Omega_{\left( \dot{x}^u, \dot{y}^u \right)} \right\}$ will be imposed in problem (P4) when $\left( x^u, y^u \right) = \left( \dot{x}^u, \dot{y}^u \right)$. Hence, $\left( \dot{x}^u, \dot{y}^u, \dot{x}^{l0}, \dot{y}^{l0} \right)$ is also feasible to problem (P4). Now let $\left( \dot{x}^u, \dot{y}^u, \dot{x}^{l0}, \dot{y}^{l0} \right)$ be a feasible solution to (P4). Following the reverse order from (52) to (50), we can show that $\left( \dot{x}^u, \dot{y}^u, \dot{x}^{l0}, \dot{y}^{l0} \right)$ is also feasible to (P2) because KKT conditions are both necessary and sufficient for optimality for LPs. We have shown above that any feasible solution to (P2) is feasible to (P4), and vice versa. In addition, the objectives of (P2) and (P4) are the same. Hence, problems (P2) and (P4) are equivalent to each other. □

We mention that, instead of utilizing $\text{Proj}_{\left( x^u, y^u \right)} P \left( y^{l,j} \right)$, i.e., the actual feasible set for the upper level when $y^l = y^{l,j}$, a different strategy that penalizes lower level constraint violations caused by any infeasible $\left( x^u, y^u, y^{l,j} \right)$ is proposed in [23]. Through using Big-M penalty coefficients, such infeasible solution will incur a large objective function value, which could make it ineligible for an optimal solution. Nevertheless, we note that selecting an appropriate Big-M is very subjective. Also, a very large M, which is desired, may seriously affect the computational performance. On the other hand, our projection-based reformulation (P4), together with the implementation method in Section 6.2, provides a rather analytical approach when the relatively complete response property is missing. Hence, it is a novel feature to the reformulation and decomposition scheme.

## 5. Algorithm

A drawback of the single-level formulation (P4) is the size of set $Y^L$ grows exponentially as the number and bound of lower-level integer variables increase. Consequently, the number of constraints (49) could be intractably large and applications to large-scale MIBLP problems may be limited. In [23], a decomposition approach based column-and-constraint generation method is developed that tries to ameliorate this issue via partial enumeration [51]. In this section, we extend that decomposition algorithm in the projection formulation context by modifying and updating the master and subproblems according to the structure of (P4), or equivalently (P0). Specifically, the algorithm always deals with a partial enumeration, where only a subset of $Y^L$ (denoted as $\underline{Y}_k^L$) is considered in each iteration $k$. Starting from an empty set $\underline{Y}_0^L = \varnothing$, we expand this subset by adding a new $y^{l,j} \in Y^L$ at the end of each iteration. This decomposition algorithm involves one master problem and two subproblems, which are presented below. As argued in [23], the master problem provides a lower bound, and the two subproblems select important $y^{l,j} \in Y^L$ that help the algorithm converge as early as possible.

### 5.1. Master problem

Let $\Theta_k^*$ be the optimal objective value of the master problem in iteration $k$. The master problem (P5) is formulated as

$$\textbf{(P5)} \quad \Theta_k^* = \min_{\substack{x^u, y^u, x^{l0}, y^{l0} \\ x^{l,j}, \pi^{l,j}}} \quad c_R^t x^u + c_Z^t y^u + d_R^t x^{l0} + d_Z^t y^{l0} \tag{53}$$

$$\text{s.t.} \quad (12), (13), \text{ and } (15)$$

$$\left[ \left( x^u, y^u \right) \in \text{Proj}_{\left( x^u, y^u \right)} P\left( y^{l,j} \right) \right]$$

$$\Rightarrow \begin{bmatrix} w_R^t x^{l0} + w_Z^t y^{l0} \geq w_R^t x^{l,j} + w_Z^t y^{l,j} \\ P_R x^{l,j} \leq s - Q_R x^u - Q_Z y^u - P_Z y^{l,j} \\ P_R^t \pi^j \geq w_R^t, x^{l,j} \perp \left( P_R^t \pi^j - w_R^t \right) \\ \pi^j \perp \left( s - Q_R x^u - Q_Z y^u - P_R x^{l,j} - P_Z y^{l,j} \right) \\ x^{l,j} \in \mathbb{R}_+^{n_R}, \pi^j \in \mathbb{R}_+^{n_L} \end{bmatrix}, \forall y^{l,j} \in \underline{Y}_k^L \subseteq Y^L \tag{54}$$

**Proposition 1** *At any given iteration $k$, $\Theta_k^*$ provides a lower bound to problem* (P4).

*Proof* Because $\underline{Y}_k^L \subseteq Y^L$, problem (P5) is a relaxation of problem (P4). Since this is a minimization problem, $\Theta_k^*$ is a valid lower bound to problem (P4). $\square$

## 5.2. Subproblem 1

Let $\left( x_k^{u,*}, y_k^{u,*} \right)$ be the optimal solution of (P5) in iteration $k$, we follow [23] and employ problem (P6) to find an optimal solution to the lower-level program at $\left( x_k^{u,*}, y_k^{u,*} \right)$.

$$\textbf{(P6)} \quad \theta_k \left( x_k^{u,*}, y_k^{u,*} \right) = \max_{x^l, y^l} \ w_R^l x^l + w_Z^l y^l \tag{55}$$

$$\text{s.t.} \ P_R x^l + P_Z y^l \le s - Q_R x_k^{u,*} - Q_Z y_k^{u,*} \tag{56}$$

$$x^l \in \mathbb{R}_+^{n_R}, y^l \in \mathbb{Z}_+^{n_Z} \tag{57}$$

where $\theta_k \left( x_k^{u,*}, y_k^{u,*} \right)$ denotes the optimal lower-level objective value at $\left( x_k^{u,*}, y_k^{u,*} \right)$. Let $\left( \hat{x}_k^l, \hat{y}_k^l \right)$ denote the optimal solution to (P6) in iteration $k$.

## 5.3. Subproblem 2

As mentioned in [23], in case where multiple optimal solutions to the lower-level program exist for a given $\left( x_k^{u,*}, y_k^{u,*} \right)$, the first subproblem (P6) may not provide the one that is desired by the upper-level objective function. Moreover, due to the existence of upper-level connecting constraint (2), $\left( \hat{x}_k^l, \hat{y}_k^l \right)$ may not even be feasible. Therefore, we modify the second subproblem presented in [23] to (P7), which either produces an optimal solution or detects an infeasible situation with $\left( x_k^{u,*}, y_k^{u,*} \right)$. It is worth pointing out that the proposition of (P7) follows a similar idea as the upper bounding program in [13] to provide an upper bound. However, we note that (P7) is a Pareto version of the lower level problem that guarantees to be optimistic, which is well-defined and provides an upper bound. Meanwhile, the upper bounding problem in [13] could be a relaxed problem due to the introduced epsilon.

$$\textbf{(P7)} \quad \Theta_{o,k} \left( x_k^{u,*}, y_k^{u,*} \right) = \min_{x^l, y^l} \ d_R^l x^l + d_Z^l y^l \tag{58}$$

$$\text{s.t.} \ (56) \text{ and } (57)$$

$$B_R x^l + B_Z y^l \le r - A_R x_k^{u,*} - A_Z y_k^{u,*} \tag{59}$$

$$w_R^l x^l + w_Z^l y^l \ge \theta_k \left( x_k^{u,*}, y_k^{u,*} \right) \tag{60}$$

where $\Theta_{o,k}\left(x_k^{u,*}, y_k^{u,*}\right)$ denotes the optimal objective value of problem (P7) if it is feasible. It is easy to see that $c_R^t x_k^{u,*} + c_Z^t y_k^{u,*} + \Theta_{o,k}\left(x_k^{u,*}, y_k^{u,*}\right)$ provides an upper bound to problem (P4).

### 5.4. Decomposition algorithm

Based on the master problem (P5) and the two subproblems (P6) and (P7), the proposed decomposition algorithm, which implements the column-and-constraint generation (CCG) method [52], is summarized below.

---

**Algorithm.** Projection-based reformulation and decomposition through CCG method

---

1  **Step 1 (Initialization)**
2      Set $LB = -\infty$, $UB = +\infty$, $\xi = 0$, $k = 0$, and $\underline{Y}_0^L \leftarrow \varnothing$.
3  **Step 2 (Lower Bounding)**
4      Solve problem (P5).
5      Denote the optimal solution as $\left(x_k^{u,*}, y_k^{u,*}, x_k^{l0,*}, y_k^{l0,*}\right)$.
6      Set $LB$ to the optimal objective value $\Theta_k^*$.
7  **Step 3 (Termination)**
8      **if** $UB - LB < \xi$, **then** Terminate and return optimal solution.
9  **Step 4 (Subproblem 1)**
10      Solve problem (P6) at $\left(x_k^{u,*}, y_k^{u,*}\right)$.
11      Denote the optimal solution as $\left(\hat{x}_k^l, \hat{y}_k^l\right)$ and optimal objective value as $\theta_k\left(x_k^{u,*}, y_k^{u,*}\right)$.
12  **Step 5 (Subproblem 2)**
13      Solve problem (P7) at $\left(x_k^{u,*}, y_k^{u,*}\right)$ and $\theta_k\left(x_k^{u,*}, y_k^{u,*}\right)$
14      **if** Feasible **then**
15          Denote the optimal solution as $\left(x_k^{l,*}, y_k^{l,*}\right)$
16          Set $UB = \min\left\{UB, c_R^t x_k^{u,*} + c_Z^t y_k^{u,*} + \Theta_{o,k}\left(x_k^{u,*}, y_k^{u,*}\right)\right\}$.
17          Set $\hat{y}_k^l = y_k^{l,*}$.
18      **else** (**Infeasible Problem**)
19          Set $\hat{y}_k^l = \hat{y}_k^l$.
20      **end**
21  **Step 6 (Tightening the Master Problem)**
22      Create new variables $\left(x^{l,j}, \pi^j\right)$ and constraint (54) corresponding to $y^{l,j} = \hat{y}_k^l$.
23      Set $\underline{Y}_{k+1}^L = \underline{Y}_{k+1}^L \cup \left\{\hat{y}_k^l\right\}$ and $k = k + 1$.
24  **Step 7 (Loop)**
25      **if** $UB - LB < \xi$, **then**

| 26 | Terminate and return the optimal solution. |
|----|---------------------------------------------|
| 27 | **else** |
| 28 | Go to step 2. |
| 29 | **end** |

**Remark 1** Assumption 1 ensures that MIBLP (P0) is feasible. Because the master problem (P5) is a relaxation of problem (P0) in all iterations, (P5) is guaranteed to be feasible.

**Remark 2** The first subproblem (P6) is always feasible. We know that at any given iteration $k$, the optimal solution to (P5) $\left( x_k^{l0,*}, y_k^{l0,*} \right)$ will be feasible to (P6).

**Remark 3** The decomposition algorithm through column-and-constraint generation provides a series of non-decreasing lower bounds. For any given iteration $k$, $\underline{Y}_k^L \subseteq \underline{Y}_{k+1}^L$. Hence, the master problem (P5) in iteration $k$ is a relaxation of that in iteration $k+1$.

## 5.5. Convergence

In this subsection, we generalize the convergence proof in [23] to our new algorithm development, noting that we consider a more general MIBLP (i.e., (P0)) and employ modified master and subproblems.

**Theorem 2** *If an MIBLP satisfies Assumption 1, $\xi = 0$ and $Y^L$ is finite, then the presented algorithm converges to the global optimal objective value of MIBLP problem* (P0) *within* $\left| Y^L \right|$ *iterations.*

*Proof* It is sufficient to show that if none of the stopping criteria are met, a new $y^{l,j} \in Y^L$ would be generated in each iteration. This is equivalent to showing that a repeated $y^{l,j} \in Y^L$ leads to either $LB = UB$ or an infeasible master problem (P5). Assume that the current iteration index is $k = l_1$ and $\left( x_{l_1}^{u,*}, y_{l_1}^{u,*}, x_{l_1}^{l0,*}, y_{l_1}^{l0,*} \right)$ is obtained in step 2 with $UB - LB > 0$. From Remark 2 we know that the first subproblem is feasible. However, the second subproblem could be infeasible.

We first consider the case where the second subproblem (P7) is feasible, so that $y_{l_1}^{l,*}$ is obtained in step 5. The proof for this case is similar to that in [23]. We further assume that $y_{l_1}^{l,*}$ was also derived in a previous iteration $k = l_0 \left( < l_1 \right)$ and included in $\underline{Y}_{l_0+1}^L$. Because $UB - LB > 0$ in iteration $l_1$, new variables $\left( x^{l,j}, \pi^j \right)$ and constraints (54) corresponding to $y^{l,j} = y_{l_1}^{l,*}$ will be added to the master problem (P5) in iteration $l_1+1$ (Note that $\hat{y}_{l_1}^l = y_{l_1}^{l,*}$). Since these variables and constraints are the same as those created and included in step 6 in iteration $l_0$, step 6 in iteration $l_1$ essentially does not change problem (P5). Consequently, it yields the same optimal values in iteration $l_1+1$ as that of iteration $l_1$, i.e., $\left(x_{l_1}^{u,*}, y_{l_1}^{u,*}, x_{l_1}^{l0,*}, y_{l_1}^{l0,*}\right) = \left(x_{l_1+1}^{u,*}, y_{l_1+1}^{u,*}, x_{l_1+1}^{l0,*}, y_{l_1+1}^{l0,*}\right)$. Hence, $LB$ does not change when the algorithm proceeds from iteration $l_1$ to $l_1+1$. In the following, we show that $LB \geq UB$ in iteration $l_1+1$.

$$LB = c_R^t x_{l_1+1}^{u,*} + c_Z^t y_{l_1+1}^{u,*} + d_R^t x_{l_1+1}^{l0,*} + d_Z^t y_{l_1+1}^{l0,*} \tag{61}$$

$$\geq c_R^t x_{l_1+1}^{u,*} + c_Z^t y_{l_1+1}^{u,*} + \min_{x^{l0}, y^{l0}} d_R^t x^{l0} + d_Z^t y^{l0} \tag{62}$$

$$\text{s.t. } B_R x^{l0} + B_Z y^{l0} \leq r - A_R x_{l_1+1}^{u,*} - A_Z y_{l_1+1}^{u,*} \tag{63}$$

$$P_R x^{l0} + P_Z y^{l0} \leq s - Q_R x_{l_1+1}^{u,*} - Q_Z y_{l_1+1}^{u,*} \tag{64}$$

$$\begin{bmatrix} w_R^t x^{l0} + w_Z^t y^{l0} \geq w_R^t x^{l,j} + w_Z^t y^{l,j} \\ P_R x^{l,j} \leq s - Q_R x_{l_1+1}^{u,*} - Q_Z y_{l_1+1}^{u,*} - P_Z y^{l,j} \\ P_R^t \pi^j \geq w_R^t, x^{l,j} \perp \left(P_R^t \pi^j - w_R^t\right) \\ \pi^j \perp \left(s - Q_R x_{l_1+1}^{u,*} - Q_Z y_{l_1+1}^{u,*} - P_R x^{l,j} - P_Z y^{l,j}\right) \\ x^{l,j} \in \mathbb{R}_+^{n_R}, \pi^j \in \mathbb{R}_+^{n_L} \end{bmatrix} \text{for } y^{l,j} = y_{l_1}^{l,*} \tag{65}$$

$$\geq c_R^t x_{l_1}^{u,*} + c_Z^t y_{l_1}^{u,*} + \min_{x^{l0}, y^{l0}} d_R^t x^{l0} + d_Z^t y^{l0} \tag{66}$$

$$\text{s.t. } B_R x^{l0} + B_Z y^{l0} \leq r - A_R x_{l_1}^{u,*} - A_Z y_{l_1}^{u,*} \tag{67}$$

$$P_R x^{l0} + P_Z y^{l0} \leq s - Q_R x_{l_1}^{u,*} - Q_Z y_{l_1}^{u,*} \tag{68}$$

$$w_R^t x^{l0} + w_Z^t y^{l0} \geq \theta_{l_1}\left(x_{l_1}^{u,*}, y_{l_1}^{u,*}\right) \tag{69}$$

$$= c_R^t x_{l_1}^{u,*} + c_Z^t y_{l_1}^{u,*} + \Theta_{o,l_1}\left(x_{l_1}^{u,*}, y_{l_1}^{u,*}\right) \tag{70}$$

$$\geq UB \tag{71}$$

Inequality (62) follows from $\left(x_{l_1+1}^{u,*}, y_{l_1+1}^{u,*}\right) \in P\left(y_{l_1}^{l,*}\right)$ and $y_{l_1}^{l,*} \in \underline{Y}_{l_1+1}^L$. Inequality (66) follows from $y_{l_1}^{l,*}$ is optimal to $\theta_{l_1}\left(x_{l_1}^{u,*}, y_{l_1}^{u,*}\right)$, and $w_R^t x^{l,j} + w_Z^t y^{l,j} = \theta_{l_1}\left(x_{l_1}^{u,*}, y_{l_1}^{u,*}\right)$ for $y^{l,j} = y_{l_1}^{l,*}$ due to the KKT conditions. Inequality (71) follows from $UB = \min\left\{UB, c_R^t x_{l_1}^{u,*} + c_Z^t y_{l_1}^{u,*} + \Theta_{o,l_1}\left(x_{l_1}^{u,*}, y_{l_1}^{u,*}\right)\right\}$. Consequently, we have $LB \geq UB$ in iteration $l_1+1$, which terminates the algorithm.

Now we consider the case where the second subproblem is infeasible, so that $\hat{y}_{l_1}^l$ is obtained in step 4 of iteration $l_1$ and $\tilde{y}_{l_1}^l = \hat{y}_{l_1}^l$. We further assume that $\hat{y}_{l_1}^l$ was also derived in a previous iteration $k = l_0 (< l_1)$ and included in $\underline{Y}^L_{l_0+1}$. In the following, we show that master problem (P5) in iteration $l_1 + 1$ is infeasible. As aforementioned, step 6 in iteration $l_1$ essentially does not change problem (P5). If we assume that the master problem (P5) in iteration $l_1 + 1$ is feasible, then we have $\left( x^{u,*}_{l_1}, y^{u,*}_{l_1}, x^{l0,*}_{l_1}, y^{l0,*}_{l_1} \right) = \left( x^{u,*}_{l_1+1}, y^{u,*}_{l_1+1}, x^{l0,*}_{l_1+1}, y^{l0,*}_{l_1+1} \right)$. Because $\left( x^{u,*}_{l_1+1}, y^{u,*}_{l_1+1} \right) \in P\left( \hat{y}^l_{l_1} \right)$, the following constraint will be imposed in problem (P5) in iteration $l_1 + 1$.

$$
\left[
\begin{array}{l}
w^t_R x^{l0} + w^t_Z y^{l0} \geq w^t_R x^{l,j} + w^t_Z y^{l,j} \\
P_R x^{l,j} \leq s - Q_R x^{u,*}_{l_1+1} - Q_Z y^{u,*}_{l_1+1} - P_Z y^{l,j} \\
P^t_R \pi^j \geq w^t_R, x^{l,j} \perp \left( P^t_R \pi^j - w^t_R \right) \\
\pi^j \perp \left( s - Q_R x^{u,*}_{l_1+1} - Q_Z y^{u,*}_{l_1+1} - P_R x^{l,j} - P_Z y^{l,j} \right) \\
x^{l,j} \in \mathbb{R}^{n_R}_+, \pi^j \in \mathbb{R}^{n_L}_+
\end{array}
\right] \text{ for } y^{l,j} = \hat{y}^l_{l_1} \quad (72)
$$

Because $\hat{y}^l_{l_1}$ is optimal to $\theta_{l_1}\left( x^{u,*}_{l_1}, y^{u,*}_{l_1} \right)$, and $w^t_R x^{l,j} + w^t_Z y^{l,j} = \theta_{l_1}\left( x^{u,*}_{l_1}, y^{u,*}_{l_1} \right)$ for $y^{l,j} = \hat{y}^l_{l_1}$ due to the KKT conditions. Hence, constraint (72) is equivalent to (69). However, we know from step 5 of iteration $l_1$ that there does not exist a tuple $\left( x^{u,*}_{l_1+1}, y^{u,*}_{l_1+1}, x^{l0,*}_{l_1+1}, y^{l0,*}_{l_1+1} \right)$ that simultaneously satisfy constraints (63), (64) and (69), which are imposed in (P5) in iteration $l_1 + 1$. This is contradictory to that the same problem is solved in iteration $l_1$. Hence, the master problem (P5) in iteration $l_1 + 1$ is infeasible. In consideration of the two cases above, we know a new $y^{l,j} \in Y^L$ would be generated in each iteration, if none of the stopping criteria are met.

Similar to that in [23], since a repeated $y^{l,j} \in Y^L$ leads to convergence of the algorithm and the fact that set $Y^L$ is finite, the algorithm will converge within $\left| Y^L \right|$ iterations. $\square$

# 6. Implementation

## 6.1. KKT-condition-based tightening constraints

As suggested by [23,13], the master problem (P5) can be tightened by introducing the following KKT-conditions related to the lower-level program at given $\left( x^u, y^u, y^{l0} \right)$. Hence, we have **(P8)** $\quad \min\limits_{\substack{x^u,y^u,x^{l0},y^{l0} \\ x^{l,j},\pi^{l,j},\tilde{x}^l,\tilde{\pi}}} \quad c_R^t x^u + c_Z^t y^u + d_R^t x^{l0} + d_Z^t y^{l0}$ (73)

$\quad$ s.t. $\quad$ (12), (13), (15), and (54)

$$w_R^t x^{l0} \geq w_R^t \tilde{x}^l \tag{74}$$

$$P_R \tilde{x}^l \leq s - Q_R x^u - Q_Z y^u - P_Z y^{l0}, P_R^t \tilde{\pi} \geq w_R^t \tag{75}$$

$$\tilde{x}^l \perp \left( P_R^t \tilde{\pi} - w_R^t \right), \tilde{\pi} \perp \left( s - Q_R x^u - Q_Z y^u - P_R \tilde{x}^l - P_Z y^{l0} \right) \tag{76}$$

$$\tilde{x}^l \in \mathbb{R}_+^{n_R}, \tilde{\pi} \in \mathbb{R}_+^{n_L} \tag{77}$$

Given that $x^u$, $y^u$ and $y^{l0}$ are variables, (P8) provides lower bound information that is parametric not only to $x^u$, $y^u$ but also to $y^{l0}$. The value of $w_R^t \tilde{x}^l + w_Z^t y^{l0}$ provides a valid lower bound support to $w_R^t x^{l0} + w_Z^t y^{l0}$, which might not be available from any fixed $y^{l,j} \in Y^L$. As mentioned in [23,13], the KKT-condition-based tightening constraints (74) – (77) help reduce the number of iterations and computational time. Therefore, unless otherwise stated, we consider KKT-condition-based tightening constraints in the decomposition algorithm for all numerical studies in this work.

### 6.2. Projection and indicator constraint

To practically handle our new projection-based formulation, we present an alternative representation for $\left( x^u, y^u \right) \in \text{Proj}_{\left( x^u, y^u \right)} P\left( y^{l,j} \right)$ in (54) in this subsection. Specifically, we use the following LP to check whether $\left( x^u, y^u \right) \in \text{Proj}_{\left( x^u, y^u \right)} P\left( y^{l,j} \right)$ or not.

For any given $y^{l,j} \in Y^L$

**(P9)** $\quad \min\limits_{x^{l,j}, t^j} \quad e^t t^j$ (78)

$\quad$ s.t. $\quad P_R x^{l,j} - t^j \leq s - Q_R x^u - Q_Z y^u - P_Z y^{l,j}$ (79)

$$x^{l,j} \in \mathbb{R}_+^{n_R}, t^j \in \mathbb{R}_+^{n_L} \tag{80}$$

where $e$ is a vector with all elements equal to 1; $t^j$ represents the relaxation variables for each constraint. The dimension of vector $t^j$ is equal to the number of lower-level constraints.

**Remark 4** $\quad$ Given $\left( x^u, y^u \right) \in \text{Proj}_{\left( x^u, y^u \right)} \Omega$ and $y^{l,j} \in Y^L$, if the optimal value of problem (P9) $e^t t^{j,*} = 0$, then $\left( x^u, y^u \right) \in \text{Proj}_{\left( x^u, y^u \right)} P\left( y^{l,j} \right)$; otherwise, $\left( x^u, y^u \right) \notin \text{Proj}_{\left( x^u, y^u \right)} P\left( y^{l,j} \right)$.

Noting that LP (P9) can be equivalently replaced by its corresponding KKT-conditions, we can replace (54) in problem (P8) with the following constraints based on Remark 4. Thus, we have

$$\textbf{(P9)} \quad \min_{\substack{x^u,y^u,x^{l0},y^{l0} \\ x^{l,j},\pi^{l,j},\bar{x}^l,\bar{\pi},t^j,\lambda^j}} \quad c_R^t x^u + c_Z^t y^u + d_R^t x^{l0} + d_Z^t y^{l0} \tag{81}$$

s.t. (12), (13), (15), (54), (74) – (77), (79) and (80)

$$\left[ e^t t^j = 0 \right] \Rightarrow \begin{bmatrix} w_R^t x^{l0} + w_Z^t y^{l0} \geq w_R^t x^{l,j} + w_Z^t y^{l,j} \\ P_R x^{l,j} \leq s - Q_R x^u - Q_Z y^u - P_Z y^{l,j} \\ P_R^t \pi^j \geq w_R^t, x^{l,j} \perp \left( P_R^t \pi^j - w_R^t \right) \\ \pi^j \perp \left( s - Q_R x^u - Q_Z y^u - P_R x^{l,j} - P_Z y^{l,j} \right) \\ x^{l,j} \in \mathbb{R}_+^{n_R}, \pi^j \in \mathbb{R}_+^{n_L} \end{bmatrix}, \forall y^{l,j} \in \underline{Y}_k^L \subseteq Y^L \tag{82}$$

$$P_R \lambda^j \geq 0, x^{l,j} \perp P_R \lambda^j, \forall y^{l,j} \in \underline{Y}_k^L \subseteq Y^L \tag{83}$$

$$e - \lambda^j \geq 0, t^j \perp \left( e - \lambda^j \right), \forall y^{l,j} \in \underline{Y}_k^L \subseteq Y^L \tag{84}$$

$$\lambda^j \perp \left( s - Q_R x^u - Q_Z y^u - P_Z y^{l,j} - P_R x^{l,j} + t^j \right), \forall y^{l,j} \in \underline{Y}_k^L \subseteq Y^L \tag{85}$$

### 6.3. Approximations

#### 6.3.1. Indicator constraint

Constraint (82) cannot be computed by off-the-shelf solvers directly. In this work, we take advantage of the special language feature – indicator constraints [1] in GAMS 24.4 [53]. We introduce a binary variable $\psi^j$ to denote whether $e^t t^j = 0$ or not. Hence, we have

$$\textbf{(P10)} \quad \min_{\substack{x^u,y^u,x^{l0},y^{l0} \\ x^{l,j},\pi^{l,j},\bar{x}^l,\bar{\pi} \\ t^j,\lambda^j,\psi^j}} \quad c_R^t x^u + c_Z^t y^u + d_R^t x^{l0} + d_Z^t y^{l0} \tag{86}$$

s.t. (12), (13), (15), (54), (74) – (77), (79), (80), (83) – (85)

$$\psi^j \Rightarrow \begin{bmatrix} w_R^t x^{l0} + w_Z^t y^{l0} \geq w_R^t x^{l,j} + w_Z^t y^{l,j} \\ P_R x^{l,j} \leq s - Q_R x^u - Q_Z y^u - P_Z y^{l,j} \\ P_R^t \pi^j \geq w_R^t, x^{l,j} \perp \left( P_R^t \pi^j - w_R^t \right) \\ \pi^j \perp \left( s - Q_R x^u - Q_Z y^u - P_R x^{l,j} - P_Z y^{l,j} \right) \\ x^{l,j} \in \mathbb{R}_+^{n_R}, \pi^j \in \mathbb{R}_+^{n_L} \end{bmatrix}, \forall y^{l,j} \in \underline{Y}_k^L \subseteq Y^L \tag{87}$$

$$\varepsilon \left( 1 - \psi^j \right) \leq e^t t^j, \psi^j \in \{0,1\}, \quad \forall y^{l,j} \in \underline{Y}_k^L \subseteq Y^L \tag{88}$$

---

[1] An indicator constraint is a way of expressing relationships among variables by specifying a binary variable to control whether or not a constraint takes effect.

where $\varepsilon$ is a very small positive number (e.g., $10^{-4}$). Constraint (88) indicates that if $e^t t^{j,*} = 0$ then $\psi^j$ is forced to 1. If $e^t t^{j,*} \geq \varepsilon$, at optimality $\psi^j$ will equal 0 as less constraints will be imposed. Note that the case $0 < e^T t^j < \varepsilon$ is excluded. Hence, formulation (P10) is an approximation of (P9).

### 6.3.2. Linearization of complementary constraints

In this work, we linearize all complementary constraints in KKT conditions by using the big-M formulation and introducing a binary variable for each complementary constraint [1]. For example,

$$0 \leq f \perp g \geq 0 \Leftrightarrow \begin{cases} 0 \leq f \leq M\delta \\ 0 \leq g \leq M(1-\delta) \\ \delta \in \{0,1\} \end{cases}$$ (89)

where $f$ and $g$ are two arbitrary equations; $M$ is a large positive number; and $\delta$ is the binary variable for complementary constraint $0 \leq f \perp g \geq 0$. It is noted that there are other approaches in handling such complementarity constraints [54-58]. We choose the big-M formulation because applying it to (P10) results in a single-level MILP with indicator constraints, which can be handled by CPLEX 12.

Assuming that $\varepsilon$ is chosen sufficiently small and $M$ is chosen sufficiently large, the proposed algorithm will converge to the optimal solution in finite iterations. In case that the infimum may not be attainable [20], the proposed algorithm converges to $\varepsilon$-optimal solutions. Interested readers are referred to [23] for more details.

## 7. Computational examples

To systematically test and evaluate our algorithm development, a very comprehensive computational study is performed on three set of examples. Specifically, we employ Example 1 to verify our algorithm by computing and comparing with the only publicly accessible MIBLP library at the time of writing this paper. We further observe that obtained results cannot reflect the full features of the proposed algorithm, given that these instances do not include upper-level continuous variables and the parameters are all assumed to be integral. We then employ Example 2 to test our algorithm on MIBLPs in the general form of (P0), which cannot be computed by the original reformulation-and-decomposition method. To demonstrate the applicability of our algorithm in solving practical problems, we employ Example 3, which is a case study on hierarchical supply chain planning. Also, a few illustrative examples are provided in the Appendix for easy understanding.

All computational experiments are performed on a PC with an Intel® Core™ i5-2400 CPU @ 3.10GHz and 8.00 GB RAM. All models and solution procedures are coded in GAMS 24.4 [59]. The resulting MILP problems are solved with CPLEX 12. The indicator constraints are programmed using the GAMS/CPLEX option files [53]. The CPLEX solver options are set as follows: $epint^2 = 0$, $eprhs^3 = 10^{-8}$. The optimality tolerances for the solver are set to 0, i.e., $optcr = 0$; $optca = 0.0$. The big-M coefficients are all set to $10^4$, and the $\varepsilon$ coefficients are set to $10^{-4}$. The tolerance of the relative gap between $UB$ and $LB$ ($\xi$) is set to $10^{-3}$.

## 7.1. Example 1

In this example, we test the proposed algorithm on the small- and medium-size instances in [7]. The computational results are presented in Table 1, where $m_R$, $m_Z$, $n_R$, and $n_Z$ are the numbers of upper-level continuous variables, upper-level integer variables, lower-level continuous variables, and lower-level integer variables. The total number of variables is $n_T = m_R + m_Z + n_R + n_Z$. The numbers of upper-level constraints ($n_U$) and lower-level constraints ($n_L$) are both set to $0.2n_T$.

The instances are put in the same order as that in [7]. The smallest instances have 0 upper-level continuous variables, 10 upper-level binary variables, an average of 5 lower-level continuous variables, an average of 5 lower-level binary variables, 4 upper-level constraints, and 4 lower-level constraints. The largest instances tested have 0 upper-level continuous variables, 110 upper-level binary variables, an average of 55 lower-level continuous variables, an average of 55 lower-level binary variables, 44 upper-level constraints, and 44 lower-level constraints. The proposed algorithm obtained the same objective function values as that in [7] for all instances. As expected, the computational time of the proposed algorithm increases as the problem dimension grows. For example, (Xu_Wang _20) instances take an average of 1 second to solve, while (Xu_Wang_220) instances take about an average of 100 seconds to solve. Furthermore, we observe that the computational time varies significantly even for problems of the same scale. For example, instance (Xu_Wang_220_1) takes 3 seconds to solve, while instance (Xu_Wang_220_10) takes about 937 seconds to solve. Another observation is that the proposed algorithm usually converges in a few number of iterations. We can see that all instances are solved within 4 iterations. From Table 1, we can see that the proposed algorithm is comparable with the branch-and-bound algorithm in [7]. For certain instances the proposed algorithm is faster, e.g., Xu_Wang _120_10; but for some instances the branch-and-bound algorithm in [7] is faster, e.g., Xu_Wang_220_10. This is because different MILP subproblems are solved, and different solution approaches are used.

**Table 1** Computational performance of the proposed algorithm on instances in [7].

| Instance | $m_R$ | $m_Z$ | $n_R$ | $n_Z$ | Proposed Alg. | | Alg. in [7] |
|---|---|---|---|---|---|---|---|
| | | | | | CPUs | # Iterations | CPUs |
| Xu_Wang_20_1 | 0 | 10 | 6 | 4 | 1 | 1 | 0 |
| Xu_Wang_20_2 | 0 | 10 | 5 | 5 | 1 | 2 | 0 |
| Xu_Wang_20_3 | 0 | 10 | 6 | 4 | 1 | 2 | 0 |
| Xu_Wang_20_4 | 0 | 10 | 4 | 6 | 1 | 2 | 1 |
| Xu_Wang_20_5 | 0 | 10 | 6 | 4 | 1 | 2 | 1 |
| Xu_Wang_20_6 | 0 | 10 | 4 | 6 | 1 | 3 | 1 |
| Xu_Wang_20_7 | 0 | 10 | 5 | 5 | 1 | 4 | 1 |
| Xu_Wang_20_8 | 0 | 10 | 6 | 4 | 1 | 2 | 2 |
| Xu_Wang_20_9 | 0 | 10 | 7 | 3 | 1 | 2 | 4 |
| Xu_Wang_20_10 | 0 | 10 | 5 | 5 | 1 | 2 | 4 |
| Xu_Wang_120_1 | 0 | 60 | 31 | 29 | 10 | 3 | 1 |
| Xu_Wang_120_2 | 0 | 60 | 25 | 35 | 5 | 2 | 2 |
| Xu_Wang_120_3 | 0 | 60 | 34 | 26 | 1 | 1 | 7 |
| Xu_Wang_120_4 | 0 | 60 | 28 | 32 | 1 | 1 | 8 |
| Xu_Wang_120_5 | 0 | 60 | 27 | 33 | 26 | 2 | 25 |
| Xu_Wang_120_6 | 0 | 60 | 33 | 27 | 1 | 1 | 30 |
| Xu_Wang_120_7 | 0 | 60 | 32 | 28 | 1 | 1 | 63 |
| Xu_Wang_120_8 | 0 | 60 | 30 | 30 | 20 | 3 | 81 |
| Xu_Wang_120_9 | 0 | 60 | 34 | 26 | 6 | 3 | 85 |
| Xu_Wang_120_10 | 0 | 60 | 35 | 25 | 3 | 2 | 154 |
| Xu_Wang_220_1 | 0 | 110 | 60 | 50 | 3 | 1 | 2 |
| Xu_Wang_220_2 | 0 | 110 | 65 | 45 | 5 | 2 | 9 |
| Xu_Wang_220_3 | 0 | 110 | 55 | 55 | 4 | 1 | 14 |
| Xu_Wang_220_4 | 0 | 110 | 60 | 50 | 1 | 1 | 15 |
| Xu_Wang_220_5 | 0 | 110 | 58 | 52 | 1 | 1 | 19 |
| Xu_Wang_220_6 | 0 | 110 | 55 | 55 | 35 | 2 | 35 |
| Xu_Wang_220_7 | 0 | 110 | 49 | 61 | 6 | 2 | 62 |
| Xu_Wang_220_8 | 0 | 110 | 56 | 54 | 21 | 2 | 75 |

| Xu_Wang_220_9 | 0 | 110 | 52 | 58 | 1 | 1 | 168 |
| Xu_Wang_220_10 | 0 | 110 | 47 | 63 | 937 | 4 | 720 |

## 7.2. Example 2

In this subsection, we test the proposed algorithm on randomly generated instances. Parameters for these instances are generated as follows. The number of upper-level variables is set equal to that of lower-level variables, i.e. $m_R + m_Z = n_R + n_Z$. The numbers of continuous and integer variables are randomly determined with equal probability. The total number of variables ($n_T = m_R + m_Z + n_R + n_Z$) is set to five levels ranging from 20 to 400. The numbers of upper-level constraints ($n_U$) and lower-level constraints ($n_L$) are set to $0.2n_T$. We use the same notations as those of problem (P0). The upper bounds of $x^u$ and $x^l$ are set to 10. The integer variables $y^u$ and $y^l$ are all assumed to be binary variables. Elements of the following matrices and vectors are real numbers randomly generated following uniform distribution. $c_R$, $c_Z$, $d_R$, $d_Z$, $w_R$, and $w_Z$ are within [–50, 50]. $A_R$, $A_Z$, $B_R$, $B_Z$, $Q_R$, $Q_Z$, $P_R$, and $P_Z$ are within [0, 10]. $r$ is within [30, 130]. $s$ is within [10, 110]. For each level of $n_T$, ten random instances are generated. The detailed inputs to the GAMs code for generating computational instances are provided in the Appendix D. Note that they are in the general form of (P0) and cannot be computed by the original reformulation-and-decomposition method.

The model statistics and computational performances corresponding to the 50 instances are summarized in Table 2. The instances that have the same total number of variables are sorted in the ascending order of computational time. The smallest instances have an average of 5 upper-level continuous variables, an average of 5 upper-level binary variables, an average of 5 lower-level continuous variables, an average of 5 lower-level binary variables, 4 upper-level constraints, and 4 lower-level constraints. The largest instances have an average of 100 upper-level continuous variables, an average of 100 upper-level binary variables, an average of 100 lower-level continuous variables, an average of 100 lower-level binary variables, 80 upper-level constraints, and 80 lower-level constraints. As expected, the computational time increases rapidly as the problem dimension grows. For example, (miblp_20) instances take an average of 1 second to solve, while (miblp_400) instances take about an average of 1 hour to solve. Furthermore, we observe that the computational time varies significantly even for problems of the same scale. For example, instance (miblp_300_1) takes 2 seconds to solve, while instance (miblp_300_10) takes about 50 minutes to solve. Another observation is that the algorithm usually converges in a few number of iterations. We can see that 49 out of the 50 instances are solved within 3 iterations, demonstrating the efficiency of the KKT-condition-based inequalities. The exception is instance (miblp_20_10) which is solved in 10 iterations due to its complexity.

**Table 2** Computational performance of the proposed algorithm on random instances.

| Instance | $m_R$ | $m_Z$ | $n_R$ | $n_Z$ | CPUs | # Iterations |
|----------|-------|-------|-------|-------|------|--------------|
| miblp_20_1 | 6 | 4 | 4 | 6 | 1 | 1 |
| miblp_20_2 | 7 | 3 | 6 | 4 | 1 | 1 |
| miblp_20_3 | 3 | 7 | 5 | 5 | 1 | 1 |
| miblp_20_4 | 6 | 4 | 5 | 5 | 1 | 3 |
| miblp_20_5 | 2 | 8 | 4 | 6 | 1 | 2 |
| miblp_20_6 | 9 | 1 | 3 | 7 | 1 | 2 |
| miblp_20_7 | 6 | 4 | 4 | 6 | 1 | 3 |
| miblp_20_8 | 3 | 7 | 6 | 4 | 1 | 2 |
| miblp_20_9 | 5 | 5 | 7 | 3 | 1 | 3 |
| miblp_20_10 | 5 | 5 | 2 | 8 | 2 | 10 |
| miblp_100_1 | 25 | 25 | 20 | 30 | 1 | 1 |
| miblp_100_2 | 26 | 24 | 26 | 24 | 1 | 1 |
| miblp_100_3 | 27 | 23 | 23 | 27 | 1 | 2 |
| miblp_100_4 | 25 | 25 | 24 | 26 | 1 | 2 |
| miblp_100_5 | 21 | 29 | 27 | 23 | 1 | 2 |
| miblp_100_6 | 30 | 20 | 28 | 22 | 1 | 2 |
| miblp_100_7 | 32 | 18 | 20 | 30 | 2 | 2 |
| miblp_100_8 | 19 | 31 | 17 | 33 | 7 | 2 |
| miblp_100_9 | 21 | 29 | 25 | 25 | 9 | 3 |
| miblp_100_10 | 22 | 28 | 28 | 22 | 13 | 3 |
| miblp_200_1 | 53 | 47 | 51 | 49 | 1 | 1 |
| miblp_200_2 | 48 | 52 | 45 | 55 | 1 | 1 |
| miblp_200_3 | 48 | 52 | 46 | 54 | 7 | 2 |
| miblp_200_4 | 40 | 60 | 58 | 42 | 19 | 2 |
| miblp_200_5 | 42 | 58 | 57 | 43 | 19 | 2 |
| miblp_200_6 | 55 | 45 | 49 | 51 | 87 | 2 |
| miblp_200_7 | 53 | 47 | 51 | 49 | 243 | 2 |
| miblp_200_8 | 51 | 49 | 49 | 51 | 268 | 2 |
| miblp_200_9 | 52 | 48 | 52 | 48 | 349 | 2 |
| miblp_200_10 | 48 | 52 | 51 | 49 | 595 | 2 |
| miblp_300_1 | 73 | 77 | 83 | 67 | 2 | 1 |
| miblp_300_2 | 80 | 70 | 68 | 82 | 2 | 1 |
| miblp_300_3 | 69 | 81 | 73 | 77 | 17 | 2 |
| miblp_300_4 | 76 | 74 | 79 | 71 | 209 | 2 |
| miblp_300_5 | 74 | 76 | 77 | 73 | 264 | 2 |
| miblp_300_6 | 79 | 71 | 75 | 75 | 290 | 2 |
| miblp_300_7 | 78 | 72 | 73 | 77 | 432 | 2 |
| miblp_300_8 | 82 | 68 | 68 | 82 | 437 | 2 |
| miblp_300_9 | 75 | 75 | 79 | 71 | 1,713 | 2 |
| miblp_300_10 | 73 | 77 | 70 | 80 | 3,016 | 2 |
| miblp_400_1 | 95 | 105 | 99 | 101 | 2 | 1 |
| miblp_400_2 | 94 | 106 | 107 | 93 | 6 | 1 |
| miblp_400_3 | 97 | 103 | 100 | 100 | 75 | 2 |
| miblp_400_4 | 98 | 102 | 104 | 96 | 93 | 1 |
| miblp_400_5 | 104 | 96 | 92 | 108 | 189 | 2 |
| miblp_400_6 | 103 | 97 | 98 | 102 | 779 | 2 |
| miblp_400_7 | 111 | 89 | 95 | 105 | 896 | 2 |
| miblp_400_8 | 99 | 101 | 97 | 103 | 8,285 | 2 |
| miblp_400_9 | 104 | 96 | 106 | 94 | 14,232 | 3 |
| miblp_400_10 | 98 | 102 | 108 | 92 | 16,573 | 2 |

### 7.3. Example 3

This problem is modified from the capacitated plant selection problem by Cao and Chen [60]. While most facility selection and production planning approaches assume centralized decision making using monolithic models, the authors addressed the problem in a decentralized manufacturing environment, where the principal firm and the auxiliary plants operate independently in an organizational hierarchy. A bilevel optimization model was proposed to separate the decision making of plant selection and production planning. Changes to the original problem include: 1) a constraint on resource limitation (e.g., capital, labor, emission cap) is added to the upper-level program; 2) a continuous capacity variable is considered as an upper-level decision variable; and 3) a fixed cost for opening a certain production line is considered. The hierarchical planning problem is formulated into an MIBLP problem which involves continuous and binary variables in both the upper- and lower-level programs. A detailed mathematical model formulation can be found in Appendix E.

We test the proposed algorithm on a total of 35 instances. By varying the number of plants and products, we consider 7 levels, each including 5 cases: (6,6), (6,8), (8,8), (8,10), (10,10), (10,12), and (12,12), where the first number denotes the number of plants and the second number denotes the number of products. The parameters are randomly generated. The demands of product $j$ ($d(j)$) are uniformly generated on $5 \times U[8,12]$. The opening cost of plant $i$ ($f(i)$) is uniformly distributed on $5 \times U[20,80]$, and the opportunity cost ($p(i)$) is generated uniformly on $0.1 \times U[4,10]$. The resource quota ($q$) is determined in each instance varying from 230 to 650. The unit production cost ($w(i)$) equals the summation of $p(i)$ and a random parameter uniformly generated on $0.1 \times U[-2,2]$. The upper bound for capacity ($cu(i)$) is uniformly distributed on $50 \times U[2,9]$. The capacity consumption ratio ($a(i,j)$) is generated as a ratio of two uniformly generated parameters, given as $0.1 \times U[7,12] / U[0,1]$. The transportation cost ($r(i,j)$) is given as a summation of three uniformly generated parameters $0.1 \times U[0,5] + 0.1 \times U[0,5] + 0.01 \times U[1,3]$. The setup cost in the upper level problem ($s(i,j)$) and the lower level problem ($g(i,j)$) are given as $round\big(0.5 \times U[20,80]\big) + 2 \times U[-3,3]$ and $round\big(0.5 \times U[20,80]\big) + 2 \times U[-3,3] + 2 \times U[-2,2]$, respectively. The resource demand for producing unit product ($e(i,j)$) is given as $0.1 \times U[0,5] + 0.1 \times U[0,5] + 0.1 \times U[1,3]$. The trivial instances that can be solved in one iteration are intentionally excluded. The detailed inputs to GAMS code for generating instances are provided in the Appendix D.

The model statistics and computational performances corresponding to the 35 instances are summarized in Table 3. The instances that have the same number of plants and products are sorted in ascending order of computational time. The smallest instances have up to 6 upper-level continuous variables, 6 upper-level binary variables, 36 lower-level continuous variables, 36 lower-level binary variables, 7 upper-level constraints, and 54 lower-level constraints. The largest instances have up to 12 upper-level continuous variables, 12 upper-level binary variables, 144 lower-level continuous variables, 144 lower-level binary variables, 13 upper-level constraints, and 180 lower-level constraints. From the model statistics, we can see that the lower-level programs are more difficult than the upper-level programs. As expected, the computational time increases as the numbers of plants and products increase. For example, (hscp_6_6) instances take an average of 1 second to solve, while (hscp_12_12) instances take an average of 16 minutes. We also observe that the computational time varies significantly even for instances with the same numbers of plants and products. For example, instance (hscp_10_10_1) takes 1 second to solve, while instance (hscp_10_10_5) takes about 4 minutes. It is noted that all instances are solved within 4 iterations. Specifically, 26 instances are solved in 2 iterations, 6 instances in 3 iterations, and 3 instances in 4 iterations.

**Table 3** Computational performance of the proposed algorithm on hierarchical supply chain planning instances.

| Instance | # plants | # products | CPUs | # Iterations |
|---|---|---|---|---|
| hscp_6_6_1 | 6 | 6 | 1 | 2 |
| hscp_6_6_2 | 6 | 6 | 1 | 2 |
| hscp_6_6_3 | 6 | 6 | 2 | 2 |
| hscp_6_6_4 | 6 | 6 | 2 | 2 |
| hscp_6_6_5 | 6 | 6 | 3 | 3 |
| hscp_6_8_1 | 6 | 8 | 1 | 2 |
| hscp_6_8_2 | 6 | 8 | 1 | 2 |
| hscp_6_8_3 | 6 | 8 | 1 | 2 |
| hscp_6_8_4 | 6 | 8 | 2 | 2 |
| hscp_6_8_5 | 6 | 8 | 102 | 3 |
| hscp_8_8_1 | 8 | 8 | 1 | 2 |
| hscp_8_8_2 | 8 | 8 | 1 | 2 |
| hscp_8_8_3 | 8 | 8 | 2 | 2 |
| hscp_8_8_4 | 8 | 8 | 2 | 2 |
| hscp_8_8_5 | 8 | 8 | 2 | 2 |
| hscp_8_10_1 | 8 | 10 | 2 | 2 |
| hscp_8_10_2 | 8 | 10 | 3 | 2 |
| hscp_8_10_3 | 8 | 10 | 5 | 2 |
| hscp_8_10_4 | 8 | 10 | 6 | 3 |
| hscp_8_10_5 | 8 | 10 | 39 | 4 |
| hscp_10_10_1 | 10 | 10 | 1 | 2 |
| hscp_10_10_2 | 10 | 10 | 7 | 2 |
| hscp_10_10_3 | 10 | 10 | 18 | 2 |
| hscp_10_10_4 | 10 | 10 | 22 | 3 |
| hscp_10_10_5 | 10 | 10 | 214 | 3 |
| hscp_10_12_1 | 10 | 12 | 4 | 2 |
| hscp_10_12_2 | 10 | 12 | 8 | 2 |
| hscp_10_12_3 | 10 | 12 | 9 | 2 |
| hscp_10_12_4 | 10 | 12 | 9 | 2 |
| hscp_10_12_5 | 10 | 12 | 117 | 3 |
| hscp_12_12_1 | 12 | 12 | 18 | 2 |
| hscp_12_12_2 | 12 | 12 | 36 | 2 |
| hscp_12_12_3 | 12 | 12 | 1,016 | 4 |
| hscp_12_12_4 | 12 | 12 | 1,214 | 2 |
| hscp_12_12_5 | 12 | 12 | 2,625 | 4 |

## 8. Conclusions

In this paper, an extended variant of the reformulation and decomposition algorithm was proposed and developed for solving a broad class of MIBLP problems. We assumed that the inducible region was nonempty and all variables had finite bounds, which guaranteed that an MIBLP is feasible and has an optimal solution. A novel projection-based single-level formulation was proposed, which accounts for MIBLPs that do not have the relatively complete response property. Based on this formulation, a decomposition algorithm through column-and-constraint generation was developed, which progressively generated stronger lower and upper bounds by iteratively solving master and subproblems. We also proved that the algorithm converges to global optimal solutions in finite iterations.

The computational performance of the proposed MIBLP solution algorithm has been comprehensively evaluated by three types of computational examples, including 30 literature instances, 50 randomly-generated numerical instances, and 35 hierarchical supply chain planning problems, formulated as MIBLPs. We conclude from the computational results that our algorithm can solve small to large scale MIBLP problems very efficiently in most cases. The future work should investigate the potential opportunity to further boost the computational performance especially when facing large number of lower-level integer variables. The exploration of MIBNLP will be another interesting direction for further exploration.

**Acknowledgements**     We greatly appreciate the helpful discussions with Professor Andreas Wächter at Department of Industrial Engineering and Management Sciences at Northwestern University. The paper has been greatly improved by the insightful and constructive feedback from the associate editor and three anonymous reviewers.  The authors acknowledge financial support from National Science Foundation (NSF) CAREER Award (CBET-1643244).

# Appendix A: Toy example 1

The following example is adapted from [8] and is a classical MIBLP problem. We use toy example 1 to verify the results of the proposed algorithm and demonstrate the solution procedure.

$$\min_{y^u, y^{l0}} \quad -y^u - 10y^{l0}$$

$$y^u \in \mathbb{Z}_+, \ y^{l0} \in \mathbb{Z}_+$$

$$\text{s.t.} \quad y^{l0} \in \arg\max_{y^l} \left\{ \begin{array}{l} -y^l : \\ -25y^u + 20y^l \le 30 \\ y^u + 2y^l \le 10 \\ 2y^u - y^l \le 15 \\ -2y^u - 10y^l \le -15 \\ y^l \in \mathbb{Z}_+ \end{array} \right\}$$

This problem does not have lower-level continuous variables, so the reformulation (P4) can be simplified and the KKT conditions are not required. In addition, since there is no upper-level constraint, the second subproblem (P7) is always feasible. Thus, it is guaranteed that a bilevel feasible solution can be obtained in each iteration. The solution procedure of the proposed algorithm is presented below, and a graphical illustration is shown in Fig. A1.

In iteration $l = 0$ we solve master problem (P5) to obtain $\left(y^{u,*}, y^{l0,*}\right) = (2,4)$ and $LB = -42$; given $y^{u,*} = 2$, we solve subproblems (P6) and (P7) to obtain $y_0^{l,*} = 2$ and $UB = -22$; at step 6, we add the following constraint to master problem(P5): $\left[1 \le y^u \le 6\right] \Rightarrow \left[y^{l0} \le 2\right]$.

In iteration $l = 1$ we solve master problem (P5) to obtain $\left(y^{u,*}, y^{l0,*}\right) = (6,2)$ and $LB = -26$; given $y^{u,*} = 6$, we solve subproblems (P6) and (P7) to obtain $y_1^{l,*} = 1$ and $UB = \min\{-22, -16\} = -22$; at step 6, we add the following constraint to master problem (P5): $\left[2.5 \le y^u \le 8\right] \Rightarrow \left[y^{l0} \le 1\right]$.

In iteration $l = 2$ we solve master problem (P7) to obtain $\left(y^{u,*}, y^{0,*}\right) = (2,2)$ and $LB = -22$; now we have $UB = LB$ so that the algorithm terminates in step 3.

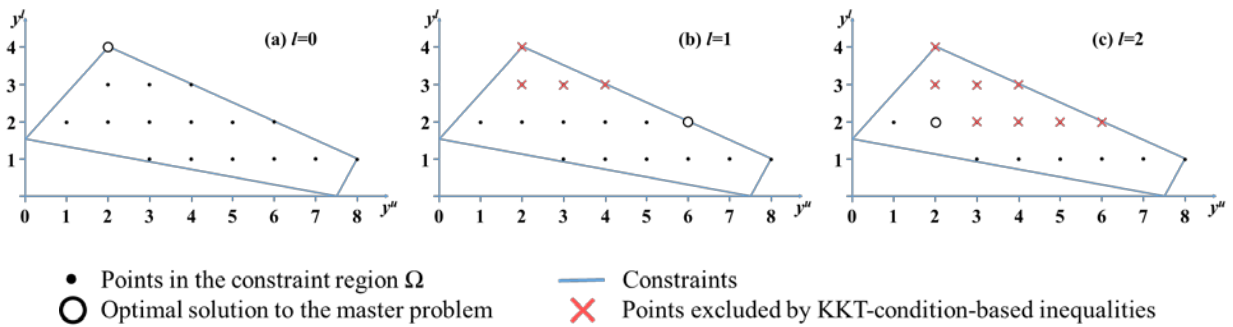

- • Points in the constraint region $\Omega$     — Constraints
- ○ Optimal solution to the master problem     × Points excluded by KKT-condition-based inequalities

**Fig. A1** The solution procedure of toy example 1.

## Appendix B: Toy example 2

The following example is adapted from [25]. We use toy example 2 to demonstrate how the proposed algorithm solves an MIBLP problem with upper-level connecting constraints. Note that this toy example and following ones in Appendices C and E cannot be computed by the original reformulation-and-decomposition method.

$$
\begin{aligned}
\min_{y^u, y^{l0}} \quad & -y^u - 2y^{l0} \\
\text{s.t.} \quad & y^u \in \mathbb{Z}_+, \, y^{l0} \in \mathbb{Z}_+ \\
& -2y^u + 3y^{l0} \le 12, \, y^u + y^{l0} \le 14 \\
& y^{l0} \in \arg\max_{y^l} \left\{ y : -3y^u + y^l \le -3, 3y^u + y^l \le 30, y^l \in \mathbb{Z}_+ \right\}
\end{aligned}
$$

This problem does not have lower-level continuous variables either, so the reformulation (P4) can be simplified and the KKT conditions are not required. However, there are two upper-level constraints that involve lower-level variables. Thus, the second subproblem (P7) could be infeasible. The solution procedure of the proposed algorithm is presented below, and a graphical illustration is shown in Fig. B1.

In iteration $l = 0$ we solve master problem (P5) to obtain $\left( y^{u,*}, y^{l0,*} \right) = \left( 6, 8 \right)$ and $LB = -22$; given $y^{u,*} = 6$, we solve the first subproblem (P6) and obtain $\hat{y}_0^l = 12$. As the second subproblem (P7) is infeasible, at step 6 we add the following constraint to master problem (P5): $\left[ 5 \le y^u \le 6 \right] \Rightarrow \left[ y^{l0} \ge 12 \right]$.

In iteration $l = 1$ we solve master problem (P5) to obtain $\left( y^{u,*}, y^{l0,*} \right) = \left( 7, 7 \right)$ and $LB = -21$; given $y^{u,*} = 7$, we solve the first subproblem (P6) and obtain $\hat{y}_1^l = 9$; but we find that the second subproblem (P7) is still infeasible; at step 6, we add the following constraint to master problem (P5): $\left[ 4 \le y^u \le 7 \right] \Rightarrow \left[ y^{l0} \ge 9 \right]$.

In iteration $l = 2$ we solve master problem (P5) to obtain $\left( y^{u,*}, y^{l0,*} \right) = \left( 8, 6 \right)$ and $LB = -20$; given $y^{u,*} = 8$, we solve subproblems (P6) and (P7) to obtain $y_2^{l,*} = 6$ and $UB = -20$; now we have $UB = LB$, so the algorithm terminates in step 7.

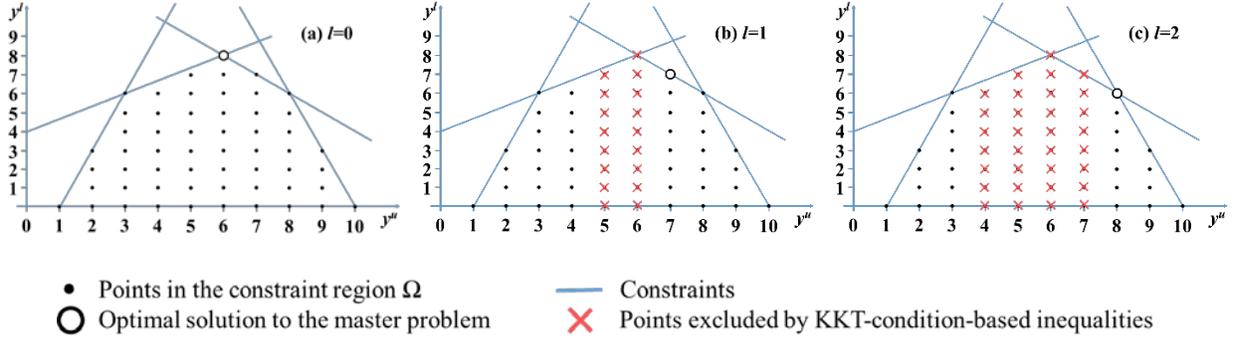

- • Points in the constraint region Ω
- ○ Optimal solution to the master problem
- — Constraints
- ✕ Points excluded by KKT-condition-based inequalities

**Fig. B1** The solution procedure of toy example 2

# Appendix C: Toy example 3

The two examples above represent a special class of MIBLPs, which include merely an upper-level integer variable and a lower-level integer variable. To show all features of the proposed algorithm while ensuring simplicity for demonstration, we propose the following illustrative example.

$$\min_{x^u, y^u, x^{l0}, y^{l0}} 20x^u - 38y^u + x^{l0} + 42y^{l0}$$

$$\text{s.t.} \quad 7y^u + 5x^{l0} + 7y^{l0} \le 62$$

$$6x^u + 9y^u + 10x^{l0} + 2y^{l0} \le 117$$

$$x^u \in \mathbb{R}_+, y^u \in \mathbb{Z}_+, x^{l0} \in \mathbb{R}_+, y^{l0} \in \mathbb{Z}_+$$

$$\left(x^{l0}, y^{l0}\right) \in \underset{x^l, y^l}{\arg\max} \ 39x^l + 27y^l$$

$$\text{s.t.} \quad 8x^u + 2x^l + 8y^l \le 53$$

$$9x^u + 2x^l + y^l \le 28$$

$$x^l \in \mathbb{R}_+, y^l \in \mathbb{Z}_+$$

This problem includes continuous and integer variables in both upper- and lower-level programs. There are two upper-level constraints involving lower-level variables. Thus, the second subproblem (P7) could be infeasible. The solution procedure of the proposed algorithm is presented below.

In iteration $l = 0$ we solve master problem (P5) to obtain $\left(x^{u,*}, y^{u,*}, x^{l0,*}, y^{l0,*}\right) = \left(2.844, 8, 1.200, 0\right)$ and $LB = -245.911$; given $\left(x^{u,*}, y^{u,*}\right) = \left(2.844, 8\right)$, we solve the first subproblem (P6) and obtain $\left(\hat{x}_0^l, \hat{y}_0^l\right) = \left(0.200, 2\right)$; the second subproblem (P7) is infeasible; at step 6 we add a set of KKT-condition-based inequalities to master problem (P5).

In iteration $l = 1$ we solve master problem (P5) to obtain $\left(x^{u,*}, y^{u,*}, x^{l0,*}, y^{l0,*}\right) = \left(2.889, 8, 1.000, 0\right)$ and $LB = -245.222$; given $\left(x^{u,*}, y^{u,*}\right) = \left(2.889, 8\right)$, we solve the first subproblem (P6) and obtain $\left(\hat{x}_1^l, \hat{y}_1^l\right) = \left(0.500, 1\right)$; we find that the second subproblem (P7) is still infeasible; at step 6, we add another set of KKT-condition-based inequalities to master problem (P5).

In iteration $l = 2$ we solve master problem (P5) to obtain $\left(x^{u,*}, y^{u,*}, x^{l0,*}, y^{l0,*}\right) = \left(3.000, 8, 0.500, 0\right)$ and $LB = -243.500$; given $\left(x^{u,*}, y^{u,*}\right) = \left(3.000, 8\right)$, we solve subproblems (P6) and (P7) to obtain $\left(x_2^{l,*}, y_2^{l,*}\right) = \left(0.500, 0\right)$ and $UB = -243.500$; now we have $UB = LB$, so the algorithm terminates in step 7.

The solution procedure above takes a total of 3 iterations. If the KKT-condition-based tightening constraints $(74) - (77)$ are not used, the algorithm takes a total of 5 iterations. Therefore, it is shown that the KKT-condition-based tightening constraints help reduce the number of iterations and computational time.

## Appendix D: Inputs for generating computational examples

The following Table 4 provides the inputs to the GAMS code for generating computational instances corresponding to example 2. We note that seed is the factor used to generate random parameters, std. stands for the standard deviation used when generating $m_R$, $m_Z$, $n_R$, and $n_Z$.

**Table 4** Inputs to the GAMS code for generating computational instances in example 2.

| Instance | seed | $0.5n_T$ | std. |
|---|---|---|---|
| miblp_20_1 | 1 | 10 | 2 |
| miblp_20_2 | 4 | 10 | 2 |
| miblp_20_3 | 1000 | 10 | 2 |
| miblp_20_4 | 7 | 10 | 2 |
| miblp_20_5 | 20 | 10 | 2 |
| miblp_20_6 | 84 | 10 | 2 |
| miblp_20_7 | 96 | 10 | 2 |
| miblp_20_8 | 5678 | 10 | 2 |
| miblp_20_9 | 79 | 10 | 2 |
| miblp_20_10 | 892 | 10 | 2 |
| miblp_100_1 | 34 | 50 | 5 |
| miblp_100_2 | 689 | 50 | 5 |
| miblp_100_3 | 1 | 50 | 5 |
| miblp_100_4 | 572 | 50 | 5 |
| miblp_100_5 | 694 | 50 | 5 |
| miblp_100_6 | 4 | 50 | 5 |
| miblp_100_7 | 42 | 50 | 5 |
| miblp_100_8 | 99 | 50 | 5 |
| miblp_100_9 | 1000 | 50 | 5 |
| miblp_100_10 | 789 | 50 | 5 |
| miblp_200_1 | 7 | 100 | 5 |
| miblp_200_2 | 377 | 100 | 5 |
| miblp_200_3 | 1065 | 100 | 5 |
| miblp_200_4 | 29 | 100 | 5 |
| miblp_200_5 | 89 | 100 | 5 |
| miblp_200_6 | 95 | 100 | 5 |
| miblp_200_7 | 232 | 100 | 5 |
| miblp_200_8 | 46 | 100 | 5 |
| miblp_200_9 | 48 | 100 | 5 |
| miblp_200_10 | 693 | 100 | 5 |
| miblp_300_1 | 10 | 150 | 5 |
| miblp_300_2 | 2 | 150 | 5 |
| miblp_300_3 | 236 | 150 | 5 |
| miblp_300_4 | 36 | 150 | 5 |
| miblp_300_5 | 25 | 150 | 5 |
| miblp_300_6 | 867 | 150 | 5 |
| miblp_300_7 | 999 | 150 | 5 |
| miblp_300_8 | 777 | 150 | 5 |
| miblp_300_9 | 239 | 150 | 5 |
| miblp_300_10 | 388 | 150 | 5 |
| miblp_400_1 | 965 | 200 | 5 |
| miblp_400_2 | 479 | 200 | 5 |
| miblp_400_3 | 374 | 200 | 5 |
| miblp_400_4 | 69 | 200 | 5 |
| miblp_400_5 | 988 | 200 | 5 |
| miblp_400_6 | 999 | 200 | 5 |
| miblp_400_7 | 111 | 200 | 5 |
| miblp_400_8 | 389 | 200 | 5 |
| miblp_400_9 | 7374 | 200 | 5 |
| miblp_400_10 | 10 | 200 | 5 |

In the following Table 5, we provide the inputs to GAMS for generating instances in example 3 from (hscp_6_6_1) through (hscp_12_12_5).

**Table 5.** Inputs to GAMS for generating instances in example 3.

| Instance | # plants | # products | seed | q |
|---|---|---|---|---|
| hscp_6_6_1 | 6 | 6 | 41257601 | 230 |
| hscp_6_6_2 | 6 | 6 | 9782 | 230 |
| hscp_6_6_3 | 6 | 6 | 18654 | 230 |
| hscp_6_6_4 | 6 | 6 | 3342 | 250 |
| hscp_6_6_5 | 6 | 6 | 22 | 260 |
| hscp_6_8_1 | 6 | 8 | 22555 | 320 |
| hscp_6_8_2 | 6 | 8 | 3611 | 350 |
| hscp_6_8_3 | 6 | 8 | 527 | 300 |
| hscp_6_8_4 | 6 | 8 | 91 | 300 |
| hscp_6_8_5 | 6 | 8 | 19123 | 360 |
| hscp_8_8_1 | 8 | 8 | 8688 | 250 |
| hscp_8_8_2 | 8 | 8 | 9651 | 300 |
| hscp_8_8_3 | 8 | 8 | 1752 | 280 |
| hscp_8_8_4 | 8 | 8 | 87422 | 250 |
| hscp_8_8_5 | 8 | 8 | 436 | 250 |
| hscp_8_10_1 | 8 | 10 | 57275355 | 400 |
| hscp_8_10_2 | 8 | 10 | 7296453 | 450 |
| hscp_8_10_3 | 8 | 10 | 72964 | 430 |
| hscp_8_10_4 | 8 | 10 | 288174 | 500 |
| hscp_8_10_5 | 8 | 10 | 2 | 450 |
| hscp_10_10_1 | 10 | 10 | 796 | 300 |
| hscp_10_10_2 | 10 | 10 | 8910 | 400 |
| hscp_10_10_3 | 10 | 10 | 23 | 350 |
| hscp_10_10_4 | 10 | 10 | 294 | 370 |
| hscp_10_10_5 | 10 | 10 | 7955 | 320 |
| hscp_10_12_1 | 10 | 12 | 89765 | 500 |
| hscp_10_12_2 | 10 | 12 | 47 | 400 |
| hscp_10_12_3 | 10 | 12 | 9364875 | 450 |
| hscp_10_12_4 | 10 | 12 | 76563 | 500 |
| hscp_10_12_5 | 10 | 12 | 3254336 | 400 |
| hscp_12_12_1 | 12 | 12 | 818 | 650 |
| hscp_12_12_2 | 12 | 12 | 97 | 350 |
| hscp_12_12_3 | 12 | 12 | 2689 | 500 |
| hscp_12_12_4 | 12 | 12 | 9434 | 480 |
| hscp_12_12_5 | 12 | 12 | 463 | 290 |

## Appendix E: Hierarchical supply chain planning model

In this section, we present the bilevel model formulation of the hierarchical supply chain planning problem adapted from [60]. Before the model is presented, we first give the notations used in the model.

*Parameters*

$a_{ij}$   capacity consumption ratio for processing product $j$ in plant $i$

$c_i^U$   upper bound of production capacity in plant $i$

$d_j$   customer demand of product $j$

$e_{ij}$   resource factor for processing product $j$ in plant $i$

$f_i$   opening cost for plant $i$

$g_{ij}$   fixed cost for opening production line $j$ in plant $i$

$p_i$   opportunity cost for unused production capacity of plant $i$ after it is opened

$q$   resource availability

$r_{ij}$   transportation cost for transferring product $j$ from plant $i$ to the principal firm

$s_{ij}$   fixed operation cost for processing product $j$ in plant $i$

$w_i$   cost to use production capacity in plant $i$

$n$   number of product types

*Continuous variables*

$Cap_i$   designated production capacity in plant $i$

$X_{ij}$   fraction of demand of product $j$ produced in plant $i$

*Binary variables*

$Y_i$       1 if plant $I$ is selected and opened; 0 otherwise

$Z_{ij}$       1 if production line for product $j$ in plant $i$ is used; 0 otherwise

With the above notations, the model for the hierarchical supply chain planning problem is formulated as follows.

$$\min \quad z_1 = \sum_i f_i Y_i + \sum_i \sum_{j \in IS_i} g_{ij} Z_{ij} + \sum_i p_i \left( Cap_i - \sum_{j \in IS_i} d_j a_{ij} X_{ij} \right) \tag{E.1}$$

$$\text{s.t.} \quad \sum_i \sum_{j \in IS_i} d_j e_{ij} X_{ij} \le q \tag{E.2}$$

$$Cap_i \le c_i^U \quad \forall i \tag{E.3}$$

$$Y_i \in \{0,1\}, Cap_i \in \mathbb{R}_+ \tag{E.4}$$

$$\min \quad z_2 = \sum_i w_i \left( \sum_{j \in IS_i} d_j a_{ij} X_{ij} \right) + \sum_i \sum_{j \in IS_i} \left( s_{ij} Z_{ij} + d_j r_{ij} X_{ij} \right) \tag{E.5}$$

$$\text{s.t.} \quad \sum_{i \in JS_j} X_{ij} = 1, \quad \forall j \tag{E.6}$$

$$\sum_{j \in IS_i} d_j a_{ij} X_{ij} \le Cap_i, \quad \forall i \tag{E.7}$$

$$\sum_{j \in IS_i} X_{ij} \le n Y_i, \quad \forall i \tag{E.8}$$

$$X_{ij} \le Z_{ij}, \quad \forall i, j \in IS_i \tag{E.9}$$

$$X_{ij} \in \mathbb{R}_+, Z_{ij} \in \{0,1\} \tag{E.10}$$

The principal firm's objective (E.1) is to minimize the sum of the plant opening cost, the production line opening cost, and the opportunity cost of over-setting production capacities. Constraint (E.2) enforces that the use of resources does not exceed their availabilities. Although only one type of resource is considered in this model, it can be easily extended to include multiple types of resources by adding an index for resources. Constraint (E.3) imposes a limitation on plant capacity. The lower-level objective function (E.5) is to minimize the operational costs, including the cost related to production capacity consumption, the fixed charge cost, and transportation costs for shipping products from auxiliary plants to the principal firm. Constraint (E.6) indicates that the demands must be fully satisfied. Constraint (E.7) indicates that production should not exceed capacity. Constraint (E.8) suggests that no product can be produced if the plant is not opened. Constraint (E.9) indicates that no product can be produced if the production line is not opened.

Constraints (E.4) and (E.10) are non-negative and binary constraints for upper- and lower-level decision variables. In this problem setting, the principal firm first determines which plant to open ($Y_i$) and the capacity to install ($Cap_i$). Then the auxiliary plants determine which production line to use ($Z_{ij}$) and the production level of each product ($X_{ij}$).